\documentclass[11pt]{article}
\usepackage{amsmath}
\usepackage{amssymb,amsmath}
\let\al=\alpha

\let\f=\frac
\let\p=\psi

\def\cA{{\cal A}}

\def\cD{{\cal D}}

\def\cF{{\cal F}}

\def\no{\noindent}
\def\na{\nabla}

\def\p{\partial}

\def\eqdefa{\buildrel\hbox{\footnotesize def}\over =}

\def\C{\mathop{\bf C\kern 0pt}\nolimits}
\def\DD{\mathop{\bf D\kern 0pt}\nolimits}
\def\K{\mathop{\bf K\kern 0pt}\nolimits}
\def\N{\mathop{\bf N\kern 0pt}\nolimits}
\def\Q{\mathop{\bf Q\kern 0pt}\nolimits}
\def\R{\mathop{\bf R\kern 0pt}\nolimits}
\def\T{\mathop{\bf T\kern 0pt}\nolimits}

\newcommand{\Z}{{\mathbf Z}}

\newcommand{\ef}{ \hfill $ \blacksquare $ \vskip 3mm}

\newcommand{\beq}{\begin{equation}}
\newcommand{\eeq}{\end{equation}}
\newcommand{\ben}{\begin{eqnarray}}
\newcommand{\een}{\end{eqnarray}}
\newcommand{\beno}{\begin{eqnarray*}}
\newcommand{\eeno}{\end{eqnarray*}}

\renewcommand{\theequation}{\thesection.\arabic{equation}}

\newtheorem{theorem}{Theorem}[section]

\newtheorem{lemma}[theorem]{Lemma}

\newtheorem{Theorem}{Theorem}[section]

\newtheorem{Proposition}[Theorem]{Proposition}

\begin{document}
\title{Well-posedness of  the Hele-Shaw-Cahn-Hilliard system}
\author{Xiaoming Wang $^\dag$ and Zhifei Zhang$^{\ddag}$\\[2mm]
{\small $ ^\dag$ Department of Mathematics, Florida State University,
Tallahassee, FL 32306-4510
  }\\
{\small E-mail: wxm@math.fsu.edu}\\[2mm]
{\small $ ^\ddag$ School of  Mathematical Science, Peking University, Beijing 100871, China}\\
{\small E-mail: zfzhang@math.pku.edu.cn}\\
{\small Corresponding author}
}

\date{\today} 
\maketitle

\begin{abstract}
We study the well-posedness of the Hele-Shaw-Cahn-Hilliard system
modeling binary fluid flow in porous media with arbitrary viscosity
contrast but matched density between the components. For initial
data in $H^s, s>\frac{d}{2}+1$, the existence and uniqueness of
solution in $C([0, T]; H^s)\cap L^2(0, T; H^{s+2})$ that is global
in time in the two dimensional case ($d=2$) and local in time in the
three dimensional case ($d=3$) are established. Several blow-up
criterions in the three dimensional case are provided as well. One
of the tools that we utilized is the Littlewood-Paley theory in
order to establish certain key commutator estimates.
\end{abstract}

\renewcommand{\theequation}{\thesection.\arabic{equation}}
\setcounter{equation}{0}

\section{Introduction}

The modeling and analysis of multi-phase fluid flow is a
fascinating, challenging and important problem
\cite{JosephRenardy1993, AMW1998}. Well-known two phase fluid
examples include the coupled atmosphere-ocean dynamical system with
water and air being the two phases, as well as the system describing
displacement of oil by water in oil reservoir (usually porous media)
\cite{Bear1988}.

A common approach to two phase flow that are macroscopically
immiscible is the sharp interface approach where the two phases are
separated by a sharp interface $\Gamma(t)$. In the case of flow in
porous media, the dynamics of the system is then governed by the two
phase Hele-Shaw (Darcy) system (Muskat problem) \cite{LLG2002,
Howison2000, SaffmanTaylor1958} together with two interface boundary
conditions: (1) continuity of the normal velocity; and (2) pressure
jump proportional to the (mean) curvature. The normal velocity of
the interface is set to be the normal velocity of the fluids.
The local in time well-posedness of the sharp interface model 
with or without surface tension is known \cite{Ambrose2004,
Ambrose2007, EscherSimonett1997}. Global in time well-posedness
with surface tension\cite{EscherSimonett1998, ConstantinPugh1993}
and 2D  without surface tension \cite{SCH2004} is also known under
the assumption that the initial data is a small perturbation of a
flat interface or a sphere.
Nevertheless, the sharp interface model encounters serious
difficulty with physically important topological changes of the
interface (possibly undefined curvature), especially in terms of
pinchoff and reconnection that are important in applications
\cite{AMW1998, LLG2002}.

As an alternative approach, one could consider the so-called phase
field models (or diffuse interface models) where an order parameter
$c$ is introduced and a capillary stress tensor is used to model the
interface between the two fluids and the forces associated
\cite{AMW1998}. The sharp interface is then replaced by a thin
transition layer and hence we avoid the difficulty of discontinuity.
In this paper, we will consider phase field approach to two phase
fluid flow with matched density in a Hele-Shaw cell or porous media.
The dynamical equations are given by the following
 Hele-Shaw-Cahn-Hilliard system \cite{LLG2002, EPalffy1997}:
\ben\label{eq:HSCH}
\left\{
\begin{array}{l}
\na\cdot u=0,\\
u=-\f1{12\eta(c)}\big(\na p-\f 1{\textbf{M}}\mu\na c\big),\\
c_t+u\cdot\na c=\f1{\textbf{Pe}}\Delta \mu,\\
c(0,x)=c_0(x),
\end{array}
\right.
\een
where $u$ is the fluid velocity, $c$ is the order parameter which is related to the concentration of the fluid, the chemical potential $\mu$ depends on the order parameter $c$ and is given by
\ben\label{eq:potential}
\mu(c)=f'_0(c)-\textbf{C}\Delta c,
\een
and \textbf{Pe} is the diffusion P\'eclet number, \textbf{C} is the Cahn number, and
\textbf{M} is a Mach number. Furthermore, $\eta(c)$ is the kinematic viscosity coefficient satisfying
\ben\label{eq:viscosity coefficient}
\eta\in C^\infty(\R^1),\quad 0<\lambda\le\eta(c)\le \Lambda<\infty,
\een
 the Helmholtz free energy  $f_0(c)$ is given by the classical double well potential
\ben\label{eq:Helmholtz free energy} f_0(c)=(c^2-1)^2. \een In the
above system (\ref{eq:HSCH}), $p$ is not the physical pressure but
the combination of certain generalized Gibbs free energy and the
gravitational potential (see \cite{LLG2002} for more details). This
model can be also viewed as the Boussinesq approximation of more
general model with arbitrary viscosity and density contrast
\cite{LLG2002}.
One may formally recover the sharp interface model 
by taking appropriate limit within the Hele-Shaw-Cahn-Hilliard
system (\ref{eq:HSCH}) \cite{LLG2002}. We will assume that the fluid
occupies the two or three dimensional torus  $\T^d, d=2,3$ for
simplicity.

Besides applications in two phase flow in porous media and Hele-Shaw
cell, certain simplified versions of this  HSCH model has been also
used in tumor growth study \cite{WLFC2008}. Moreover,
unconditionally stable schemes has been developed \cite{Wise2010}
and the existence of certain type of weak solutions (without
uniqueness) is also derived \cite{FengWise2010} for the case with
matched density and viscosity.

The goal of this manuscript is to study the well-posedness of the
matched density Hele-Shaw-Cahn-Hilliard system (\ref{eq:HSCH}) with
arbitrary viscosity contrast.

The Hele-Shaw-Cahn-Hilliard system can be formally viewed as an
appropriate limit of the classical Navier-Stokes-Cahn-Hilliard
system \cite{AMW1998, LLG2002, HH1977} which is a popular phase
field model for two phase flow  although no rigorous justification
is known yet. There are a lot of works on the Navier-Stokes-Cahn-Hilliard system including local in time well-posedness in 2 and 3 dimensional  and global in time well-posedness in 2D under various assumptions 
\cite{Abels2009, Boyer1999}. In fact the global in time
well-posedness of the 2D Navier-Stokes-Cahn-Hilliard system is
recently resolved \cite{Abels2009} using a very different set of
tools than employed here. Mathematically speaking, the difficulty
associated with the Hele-Shaw-Cahn-Hilliard is about the same as
those associated with the Navier-Stokes-Cahn-Hilliard: we gain the
advantage of dropping the nonlinear advection term in the velocity
equation but also lose the regularizing viscosity term; and their
scaling behaviors are very similar. We refer to
\cite{Lin1995,Lin2002,Xu2010, AMW1998} and references
therein for more related works on the Navier-Stokes-Cahn-Hilliard system.. 

The rest of the paper is organized as follows. We prove a key estimate on the ``pressure" in the second section. This estimate is nontrivial due to the variable coefficient introduced with the mismatched viscosity. New estimates on certain commutator operators in fractional derivative spaces are needed and they are derived in the Appendix. In section three we present the local in time well-posedness based on certain modified Galerkin approximation of the HSCH system and the ``pressure" estimate from section 2.  In section 4 we provide a Beale-Kao-Majda type blow-up criterion and prove that the system is global in time well-posed in the two dimensional case. We provide a refined blow-up criterion in the 3D case in section 5.


\section{The estimate of the pressure}
\setcounter{equation}{0}

In this section, we present the estimate of the modified pressure $p$. Taking the divergence for
the second equation of (\ref{eq:HSCH}), we find that
\ben\label{eq:pressure}
\textrm{div}\Big(\f1{\eta(c)}\na p\Big)=\textrm{div}\Big(\f1{\eta(c)}\mu(c)\na c\Big)\eqdefa \textrm{div}F(c).
\een
This variable coefficient problem is dealt with utilizing commutator estimates that we derived in the Appendix. The commutator estimates themselves are derived utilizing Littlewood-Paley decomposition.

\begin{Proposition}\label{lem:pressure}
Let $s\ge 0$ and $c\in H^{s+2}(\T^d)$, and $p$ be a smooth solution of (\ref{eq:pressure}).
Then the solution $p$ satisfies
\ben\label{eq:pressure-est}
\|\na p\|_{H^s}\le \cF(\|c\|_{L^\infty})\big(1+\|\na c\|_{L^\infty}\big)\big(1+\|c\|_{H^2}\big)^{k}\|c\|_{H^{s+2}}.
\een
Here $k=[2s]+1$ and $\cF$ is an increasing function on $\R^+$.
\end{Proposition}

\no{\bf Proof.}\, Thanks to (\ref{eq:viscosity coefficient}), a straightforward energy estimate yields that
\ben\label{eq:pressure-L2}
\|\na p\|_{L^2}\le C\|\mu(c)\|_{L^2}\|\na c\|_{L^\infty}
\le C(1+\|c\|_{L^\infty}^2)\|\na c\|_{L^\infty}\|c\|_{H^2}.
\een

Taking the operator $\langle D\rangle^s$ to (\ref{eq:pressure}) to obtain
\beno
\textrm{div}\Big(\f1{\eta(c)}\na \langle D\rangle^sp\Big)&=&\textrm{div}\langle D\rangle^s\Big(\f1{\eta(c)}\mu(c)\na c\Big)
-\textrm{div}\Big(\langle D\rangle^s\Big(\f1{\eta(c)}\na p\Big)-\Big(\f1{\eta(c)}\na \langle D\rangle^sp\Big)\Big)\\
&=& \textrm{div}(A+B),
\eeno
from which and the energy estimate, we infer that
\beno
\|\na p\|_{H^s}\le C\big(\|A\|_{L^2}+\|B\|_{L^2}\big).
\eeno

Due to the definition of $\mu(c)$, we have
\beno
\f1{\eta(c)}\mu(c)\na c=\f1{\eta(c)}f'_0(c)\na c-\textbf{C}\f1{\eta(c)}\Delta c\na c
=\na g_1(c)-\Delta c\na g_2(c),
\eeno
for some $g_1, g_2$ with $g_1(0)=g_2(0)=0$. We have by Lemma \ref{lem:nonlinear} that
\beno
\|\langle D\rangle^s\na g_1(c)\|_{L^2}\le \cF(\|c\|_{L^\infty})\|c\|_{H^{s+1}},
\eeno
and using Bony's decomposition to write
\beno
\Delta c\na g_2(c)&=&T_{\Delta c}\na g_2(c)+\widetilde{R}(\Delta c, \na g_2(c))\\
&=&\textrm{div}T_{\na c}\na g_2(c)-T_{\na c}\cdot\na \na g_2(c)+\widetilde{R}(\Delta c, \na g_2(c)),
\eeno
then from the proof of Lemma \ref{lem:product}, it is easy to see that
\beno
\|\langle D\rangle^s\Delta c\na g_2(c)\|_{L^2}\le \cF(\|c\|_{L^\infty})\|\na c\|_{L^\infty}\|c\|_{H^{s+2}}.
\eeno
Thus we obtain
\beno
\|A\|_{L^2}\le \cF(\|c\|_{L^\infty})\big(1+\|\na c\|_{L^\infty}\big)\|c\|_{H^{s+2}}.
\eeno
and by Lemma \ref{lem:commutator}-\ref{lem:nonlinear} and (\ref{eq:pressure-L2}), for $s\in (0,1]$,
\beno
\|B\|_{L^2}\le \cF(\|c\|_{L^\infty})\|c\|_{H^{s+2}}\|\na p\|_{L^2}
\le \cF(\|c\|_{L^\infty})\|\na c\|_{L^\infty}\|c\|_{H^2}\|c\|_{H^{s+2}}.
\eeno
Thus we obtain that for $s\in (0,1]$,
\ben\label{eq:pressure-H1}
\|\na p\|_{H^s}\le \cF(\|c\|_{L^\infty})\big(1+\|\na c\|_{L^\infty}\big)\big(1+\|c\|_{H^2}\big)\|c\|_{H^{s+2}},
\een

For general $s$, we will prove it by the induction argument. Let us assume that for $s\in (\f {k-1}2,\f k2]$, we have
\beno
\|\na p\|_{H^s}\le \cF(\|c\|_{L^\infty})\big(1+\|\na c\|_{L^\infty}\big)\big(1+\|c\|_{H^2}\big)^k\|c\|_{H^{s+2}}.
\eeno
Note that (\ref{eq:pressure-H1}) means that the cases of $k=1,2$ hold. Now let us assume $s\in (\f k 2,\f {k+1} 2]$.
We infer from Lemma \ref{lem:commutator} and Lemma \ref{lem:nonlinear} that
\beno
\|B\|_{L^2}\le \cF(\|c\|_{L^\infty})\big(\|c\|_{H^{s+2}}\|\na p\|_{L^2}+\|c\|_{H^2}\|\na p\|_{H^{s-\f 12}}\big).
\eeno
Then from (\ref{eq:pressure-L2}) and the induction assumption, it follows that
\beno
\|B\|_{L^2}\le \cF(\|c\|_{L^\infty})\big(1+\|\na c\|_{L^\infty}\big)\big(1+\|c\|_{H^2}\big)^{k+1}\|c\|_{H^{s+2}}.
\eeno
Thus for $s\in (\f k2,\f {k+1}2]$, we have
\beno
\|\na p\|_{H^s}\le \cF(\|c\|_{L^\infty})\big(1+\|\na c\|_{L^\infty}\big)\big(1+\|c\|_{H^2}\big)^{k+1}\|c\|_{H^{s+2}}.
\eeno
This completes the proof of Lemma \ref{lem:pressure}.\ef

\medskip

\noindent{\it Remark:} Instead of relying on the estimates from the appendix which depend on the Littlewood-Paley theory, classical energy method might work as well if we are content with
less sharp and less general results. For instance,
if $\na c\in L^\infty(\T^d)$ and $c\in H^k(\T^d)$ for $k\in \Z^+$,  classical elliptic
estimates may lead to
\beno
\|\na p\|_{H^k}\le C\|F(c)\|_{H^k},
\eeno
where $C$ depends on $\|\na c\|_{L^\infty}$ and $\|c\|_{H^k}$. 
And a straightforward product estimate gives
\beno
\|F(c)\|_{H^k}\le C(\|c\|_{L^\infty})\big(\|c\|_{H^{k+1}}+\|\na c\|_{L^\infty}\|c\|_{H^{k+2}}+
\|\Delta c\|_{L^\infty}\|c\|_{H^{k+1}}\big).
\eeno
This estimate is enough to obtain the local well-posedness of the system (\ref{eq:HSCH}) and global well-posedness in the 2D case in the space of
$$c\in C([0,T];H^2(\T^d))\cap L^2(0,T;H^{4}(\T^d)),\quad u\in C([0,T];L^2(\T^d))\cap L^2(0,T;H^{2}(\T^d))$$
when combined with the $L^\infty(H^2)\cap L^2(H^4)$ a priori
estimates from (Theorem \ref{thm:BKM criterion and GWP}) for initial
data in $H^k, k>2$. However, in order to obtain the sharp blow-up
criterion which in particular implies the global existence of the
2-D system in general Sobolev spaces  as specified in Theorem
\ref{thm:local well-posedness}, we need to establish the  refined
pressure estimate (\ref{eq:pressure-est}). Notice that
(\ref{eq:pressure-est}) is established for general (Hilbert) Sobolev
spaces, and only a linear  in $\|\na c\|_{L^\infty}$ factor in the
estimate appears in contrast to pure energy estimates.

\section{Local well-posedness}
\setcounter{equation}{0}

In this section we prove the local well-posedness of the Hele-Shaw-Cahn-Hilliard system. The procedure is mostly standard except the pressure estimate.

\begin{Theorem}\label{thm:local well-posedness}
Let $c_0(x)\in H^s(\T^d)$ for $s>\f d2+1$. Then there exists $T>0$ such that
the system (\ref{eq:HSCH}) has a unique solution $(c,u)$ in $[0,T]$ with
\beno
c\in C([0,T];H^s(\T^d))\cap L^2(0,T;H^{s+2}(\T^d)),\quad u\in C([0,T];H^{s-2}(\T^d))\cap L^2(0,T;H^{s}(\T^d));
\eeno
and satisfying the following energy estimate
\ben\label{eq:energy estimate}
\|c(t)\|_{H^s}^2+\int_0^t\|c(\tau)\|_{H^{s+2}}^2d\tau
\le \|c_0\|_{H^s}\exp\big(\int_0^tG(\tau)d\tau\big).
\een
for $t\in [0,T]$, where
\beno
G(t)=\cF(\|c\|_{L^\infty})\big(1+\|\na c\|_{L^\infty}\big)^2\big(
\|\na c\|_{L^\infty}+\|c\|_{H^{3}}^{\f{d-2}2}\big)^2\big(1+\|c\|_{H^2}\big)^{2([2s]+1)}.
\eeno
\end{Theorem}

\no{\bf Proof.} We will use the energy method to prove Theorem \ref{thm:local well-posedness}.\vspace{0.1cm}

\textbf{Step 1}. Construction of an approximate solution sequence.

The construction of the approximate solutions is based on Galerkin method. Let us define the operator $P_n$ by
\beno
P_nf(x)=\sum_{|k|\le n}f_ke^{2\pi ik\cdot x}, \quad f_k=\int_{\T^d}f(x)e^{-2\pi i k\cdot x}dx.
\eeno
Then we consider the following approximate system of (\ref{eq:HSCH}):
\ben\label{eq:approximate system}
\left\{
\begin{array}{l}
\na\cdot u_n=0,\\
u_n=-\f1{12\eta(P_n c_n)}\big(\na p_n-\f 1{\textbf{M}}\mu(P_nc_n)\na P_n c_n\big),\\
\p_tc_{n}+P_n(u_n\cdot\na P_n c_n)=\f1{\textbf{Pe}}\Delta P_n\mu(P_n c_n),\\
c_n(0,x)=P_nc_0(x).
\end{array}
\right.
\een
It is easy to see that
\beno
\|\Delta P_n\mu(P_n c_n^1)-\Delta P_n\mu(P_n c_n^2)\|_{L^2}\le C(n,\|c_n^1\|_{L^2},\|c_n^2\|_{L^2} )\|c_n^1-c_n^2\|_{L^2}.
\eeno
Taking the divergence to the second equation in (\ref{eq:approximate system}) gives
\beno
\textrm{div}\Big(\f1{12\eta(P_n c_n)}\na p_n\Big)=\f 1{\textbf{M}}\textrm{div}\Big(\f1{12\eta(P_n c_n)}\mu(P_nc_n)\na P_n c_n\Big).
\eeno
Thanks to (\ref{eq:viscosity coefficient}), straightforward energy estimate yields that
\beno
\|\na p_n\|_{L^2}\le C(n,\|c_n\|_{L^2})\|c_n\|_{L^2},
\eeno
thus we infer from the second equation of (\ref{eq:approximate system}) that
\beno
\|u_n\|_{L^2}\le C(n,\|c_n\|_{L^2})\|c_n\|_{L^2}.
\eeno
Therefore, we have
\beno
\|P_n(u_n^1\cdot\na P_n c_n^1)-P_n(u_n^2\cdot\na P_n c_n^2)\|_{L^2}\le C(n,\|c_n^1\|_{L^2},\|c_n^2\|_{L^2} )\|c_n^1-c_n^2\|_{L^2}.
\eeno
Thus, the Cauchy-Lipschtiz theorem ensures that there exists $T_n>0$ such that the approximate
system (\ref{eq:approximate system}) has a unique solution $c_n\in C([0,T_n];L^2(\T^d))$. Note that
$P_n^2=P_n$, $P_nc_n$ is also a solution of (\ref{eq:approximate system}). So the uniqueness implies
that $P_nc_n=c_n$. Thus, the approximate system (\ref{eq:approximate system}) reduces to
\ben\label{eq:approximate system-new}
\left\{
\begin{array}{l}
\na\cdot u_n=0,\\
u_n=-\f1{12\eta(c_n)}\big(\na p_n-\f 1{\textbf{M}}\mu(c_n)\na c_n\big),\\
\p_tc_{n}+P_n(u_n\cdot\na c_n)=\f1{\textbf{Pe}}\Delta P_n\mu(c_n),\\
c_n(0,x)=P_nc_0(x).
\end{array}
\right.
\een
In what follows, we denote $T_n^*$ by the maximal existence time of the solution $c_n$. Due to $P_nc_n=c_n$,
the solution $c_n$ is in fact smooth.
\vspace{0.2cm}

\textbf{Step 2}. Energy estimates.

Although the HSCH model (\ref{eq:HSCH}) has a natural energy (which is somewhat equivalent to $H^1$ estimate, see \cite{LLG2002, Wise2010} and section 4 below), it is not sufficient for the strong solution. Therefore we have to derive estimates in Sobolev spaces with higher derivatives.

For this purpose we take the $H^s(\T^d)$ inner product of the third
equation (\ref{eq:approximate system-new}) with $c_n$ and obtain
\ben\label{eq:energy}
\f12\f{d}{dt}\|c_n\|_{H^s}^2-\f1{\textrm{Pe}}\big(\Delta
P_n\mu(c_n), c_n\big)_{H^s}=-\big(u_n\cdot\na c_n, c_n\big)_{H^s}.
\een Due to (\ref{eq:potential}), we see that \beno -\big(\Delta
P_n\mu(c_n), c_n\big)_{H^s}=\textrm{C}\|\Delta c_n\|^2_{H^s}
-\big(\Delta f_0'(c_n), c_n\big)_{H^s}. \eeno We deduce, thanks to
Lemma \ref{lem:product} that \ben\label{eq:f-part} \big|\big(\Delta
f_0'(c_n), c_n\big)_{H^s}\big|\le \|f_0'(c_n)\|_{H^s}\|\Delta
c_n\|_{H^s} \le
C\big(1+\|c_n\|_{L^\infty}^2\big)\|c_n\|_{H^s}\|\Delta c_n\|_{H^s}.
\een and by Lemma \ref{lem:product} with $\sigma=1$,
\ben\label{eq:u-inner product}
\big|\big(u_n\cdot\na c_n, c_n\big)_{H^s}\big|&\le& \|u_n\cdot\na c_n\|_{H^s}\|c_n\|_{H^s}\nonumber\\
&\le& C\big(\|u_n\|_{H^s}\|\na c_n\|_{L^\infty}+\|u_n\|_{H^{\f d 2-1}}\|\na c_n\|_{H^{s+1}}\big)\|c_n\|_{H^s}.
\een
Thanks to (\ref{eq:approximate system-new}), we find that
\ben\label{eq:u-Hs}
\|u_n\|_{H^s}\le C\big(\|\f1{\eta(c_n)}\na p\|_{H^s}+\|\f1{\eta(c_n)}\mu(c_n)\na c_n\|_{H^s}\big).
\een
By Lemma \ref{lem:product}, Lemma \ref{lem:nonlinear} and Proposition \ref{lem:pressure},
the first term on the right hand side of (\ref{eq:u-Hs})  is bounded by
\beno
&&\cF(\|c_n\|_{L^\infty})\big(\|c_n\|_{H^{s+\f d2}}\|\na p\|_{L^2}+\|\na p\|_{H^s}\big)\\
&& \le \cF(\|c_n\|_{L^\infty})\big(1+\|\na c_n\|_{L^\infty}\big)\big(1+\|c_n\|_{H^2}\big)^{[2s]+1}\|c_n\|_{H^{s+2}},
\eeno
and the second term is bounded by
\beno
\cF(\|c_n\|_{L^\infty})\big(1+\|\na c_n\|_{L^\infty}\big)\|c_n\|_{H^{s+2}}.
\eeno
Thus we obtain
\beno
\|u_n\|_{H^s}\le \cF(\|c_n\|_{L^\infty})\big(1+\|\na c_n\|_{L^\infty}\big)\big(1+\|c_n\|_{H^2}\big)^{[2s]+1}\|c_n\|_{H^{s+2}},
\eeno
and similarly,
\beno
\|u_n\|_{H^{\f d2-1}}\le \cF(\|c_n\|_{L^\infty})\big(1+\|\na c_n\|_{L^\infty}\big)
\big(1+\|c_n\|_{H^2}\big)^{d-1}\|c_n\|_{H^{\f d2+1}},
\eeno
from which and (\ref{eq:u-inner product}), we infer that
\ben\label{eq:u-part}
\big|\big(u_n\cdot\na c_n, c_n\big)_{H^s}\big|&\le& \cF(\|c_n\|_{L^\infty})\big(1+\|\na c_n\|_{L^\infty}\big)\nonumber\\
&&\times\big(
\|\na c_n\|_{L^\infty}+\|c_n\|_{H^{3}}^{\f{d-2}2}\big)\big(1+\|c_n\|_{H^2}\big)^{[2s]+1}\|c_n\|_{H^{s+2}}\|c_n\|_{H^s}.
\een
Here we used the following interpolation inequality:
\beno
\|c_n\|_{H^{\f d2+1}}\le \|c_n\|_{H^{2}}^{2-\f d2}\|c_n\|_{H^{3}}^{\f d2-1}.
\eeno

Plugging (\ref{eq:f-part}) and (\ref{eq:u-part}) into (\ref{eq:energy}) yields that
\beno
&&\f12\f{d}{dt}\|c_n\|_{H^s}^2+\f {\textbf{C}} {\textbf{Pe}}\|\Delta c_n\|_{H^s}^2\\
&&\le \cF(\|c_n\|_{L^\infty})\big(1+\|\na c_n\|_{L^\infty}\big)\big(
\|\na c_n\|_{L^\infty}+\|c_n\|_{H^{3}}^{\f{d-2}2}\big)
\big(1+\|c_n\|_{H^2}\big)^{[2s]+1}\|c_n\|_{H^{s+2}}\|c_n\|_{H^s},
\eeno
which along with Young's inequality implies that
\beno
&&\f{d}{dt}\|c_n\|_{H^s}^2+\|c_n\|_{H^{s+2}}^2\\
&&\le \cF(\|c_n\|_{L^\infty})\big(1+\|\na c_n\|_{L^\infty}\big)^2\big(
\|\na c_n\|_{L^\infty}+\|c_n\|_{H^{3}}^{\f{d-2}2}\big)^2\big(1+\|c_n\|_{H^2}\big)^{2([2s]+1)}\|c_n\|_{H^s}^2.
\eeno
Then Gronwall's inequality applied gives
\ben\label{eq:energy estimate-app}
E^s_n(t)\eqdefa \|c_n(t)\|_{H^s}^2+\int_0^t\|c_n(\tau)\|_{H^{s+2}}^2d\tau
\le \|c_0\|_{H^s}\exp\big(\int_0^tG_n(\tau)d\tau\big)
\een
for $t\in [0,T_n^*)$, where
\beno
G_n(t)=\cF(\|c_n\|_{L^\infty})\big(1+\|\na c_n\|_{L^\infty}\big)^2\big(
\|\na c_n\|_{L^\infty}+\|c_n\|_{H^{3}}^{\f{d-2}2}\big)^2\big(1+\|c_n\|_{H^2}\big)^{2([2s]+1)}.
\eeno

\textbf{Step 3.} Uniform estimates and existence of the solution.

Let us define
\beno
\widetilde{T}_n^*\eqdefa \sup\big\{t\in [0,T_n^*): E_n^s(\tau)\le 2\|c_0\|_{H^s}^2\textrm{ for }\tau\in [0,t]\big\}.
\eeno
From (\ref{eq:energy estimate-app}) and Sobolev embedding, we find that
\beno
E_n^s(t)&\le& \|c_0\|_{H^s}\exp\big(\cA(\|c_0\|_{H^s})\int_0^t(1+\|c(\tau)\|_{H^3}^{d-2})d\tau\big)\\
&\le& \|c_0\|_{H^s}\exp\big(\cA(\|c_0\|_{H^s})(t+t^\f12)\big),
\quad t\in [0,\widetilde{T}^*_n).
\eeno
Here $\cA(\cdot)$ is some increasing function. Take $T$ be small enough such that
\beno
\exp\big(\cA(\|c_0\|_{H^s})(T+T^\f12)\big)\le \f32.
\eeno
Now we can conclude that $\widetilde{T}_n^*\ge T$. Otherwise, we have
\beno
E_n^s(t)\le \f32\|c_0\|_{H^s}^2\quad \textrm{for} \quad t\in [0,\widetilde{T}^*_n],
\eeno
which contradicts with the definition of $\widetilde{T}_n^*$. Thus the approximate solution $(c_n,u_n)$ exists
on $[0,T]$ and satisfies the following uniform estimate
\ben\label{eq:unifrom estimate}
\|c_n(t)\|_{H^s}^2+\int_0^t\|c_n(\tau)\|_{H^{s+2}}^2d\tau\le 2\|c_0\|_{H^s}
\een
for $t\in [0,T]$. On the other hand, it is easy to verify from the third equation of (\ref{eq:approximate system-new})
that $\p_t c_n$ is uniformly bounded in $L^2(0,T;H^{s-2}(\T^d))$.
Thus, Lions-Aubin's compactness theorem ensures that there exist a subsequence $(c_{n_k}, u_{n_k})_{k}$ of $(c_n,u_n)_n$ and
a function $c\in L^\infty(0,T;H^s(\T^d))\cap L^2(0,T;H^{s+2}(\T^d))$
and $u\in L^\infty(0,T;H^{s-2}(\T^d))\cap L^2(0,T;H^{s}(\T^d))$ such that
\beno
&&c_{n_k}\longrightarrow c, \qquad \textrm{in}\quad L^2(0,T; H^{s'+2}(\T^d)),\\
&&u_{n_k}\longrightarrow u, \qquad \textrm{in}\quad L^2(0,T; H^{s'}(\T^d)),
\eeno
as $k\rightarrow +\infty$, for any $s'<s$. Then passing to limit in (\ref{eq:approximate system-new}),
it is easy to see that $(c,u)$ satisfies (\ref{eq:HSCH}) in the weak sense
and $(c,u)$ satisfies (\ref{eq:energy estimate}).\vspace{0.1cm}

{\bf Step 4.} Continuity in time of the solution.

Revisiting the proof of (\ref{eq:energy estimate-app}), we can in fact obtain better estimate for $c_n$ (thus for $c$):
\beno
\|c\|_{\widetilde{L}^\infty(0,T; H^s(\T^d))}^2\eqdefa \sum_{j\ge -1}2^{2js}\|\Delta_j c\|_{L^\infty(0,T;L^2)}^2
\le C,
\eeno
which will imply $c\in C([0,T];H^s(\T^d))$. In fact, for any $\varepsilon>0$,
take $N$ big enough such that
\beno
\sum_{j>N}2^{2js}\|\Delta_j c\|_{L^\infty(0,T;L^2)}^2\le \f \varepsilon 4.
\eeno
For any $t\in (0,T)$ and $\delta$ such that $t+\delta\in [0,T]$, we have
\beno
\|c(t+\delta)-c(t)\|_{H^s}^2&\le& \sum_{j=-1}^N2^{2js}\|\Delta_jc(t+\delta)-\Delta_jc(t)\|_{L^2}^2+\f \varepsilon 2\\
&\le& \sum_{j=-1}^N2^{2js}|\delta|\|\p_t c\|_{L^2(0,T;L^2)}^2+\f \varepsilon 2\\
&\le& 2N2^{2N}\|\p_t c\|_{L^2(0,T;L^2)}^2|\delta|+\f \varepsilon 2.
\eeno
Thus for $|\delta|$ small enough, we have
\beno
\|c(t+\delta)-c(t)\|^2\le \varepsilon.
\eeno
That is, $c(t)$ is continuous in $H^s(\T^d)$ at the time $t$, thus  so does $u$.\vspace{0.2cm}

{\bf Step 5.} Uniqueness of the solution

Assume that $(c_1,u_1)$ and $(c_2,u_2)$ are two solutions of (\ref{eq:HSCH}) with the same initial data.
We introduce the difference of two solutions:
\beno
\delta_c=c_1-c_2, \quad \delta_u=u_1-u_2.
\eeno
Then $(\delta_c, \delta_u)$ satisfies
\beno
\left\{
\begin{array}{l}
\p_t{\delta_c}+u_1\cdot\na \delta_c+\delta_u\cdot \na c_2=\f1{\textbf{Pe}}\Delta\big(\mu(c_1)-\mu(c_2)\big) ,\\
\delta_u=\f {\eta(c_1)-\eta(c_2)}{12\eta(c_1)\eta(c_2)}\big(\na p_1-\f 1{\textbf{M}}\mu(c_1)\na c_1\big)
-\f1{12\eta(c_2)}\big(\na(p_1-p_2)-\f 1{\textbf{M}}(\mu(c_1)\na c_1-\mu(c_2)\na c_2)\big),\\
\delta_c(0)=0.
\end{array}
\right.
\eeno

Making $L^2(\T^d)$ energy estimate yields that
\beno
\f12\f d{dt}\|\delta_c\|_{L^2}^2+\f {\textbf{C}} {\textbf{Pe}}\|\Delta \delta_c\|_{L^2}^2
&\le& \f 1 {\textbf{Pe}}\big(\Delta(f_0'(c_1)-f_0'(c_2),\delta_c)\big)_{L^2}-\big(\delta_u\cdot\na c_2,\delta_c\big)_{L^2}\\
&\le& C\big(\|\Delta \delta_c\|_{L^2}+\|\delta_u\|_{L^2}\big)\|\delta_c\|_{L^2}.
\eeno
On the other hand, we can deduce from the equation of $\delta_u$ that
\beno
\|\delta_u\|_{L^2}&\le& C\big(\|\delta_c\|_{L^2}+\|\na (p_1-p_2)\|_{L^2}+\|\Delta \delta_c\|_{L^2}\big)\\
&\le& C\big(\|\delta_c\|_{L^2}+\|\Delta \delta_c\|_{L^2}\big).
\eeno
Thus we obtain
\beno
\f d{dt}\|\delta_c\|_{L^2}^2\le C\|\delta_c\|_{L^2}^2,\quad \|\delta_c(0)\|=0,
\eeno
which along with Gronwall's inequality implies $\delta_c=0$, and the uniqueness follows.\ef

\section{Blow-up criterion and global existence in 2D}

In this section we prove a Beale-Kato-Majda type blow-up criterion \cite{MB2002}  for the Hele-Shaw-Cahn-Hilliard system.
As an application, we obtain the global well-posedness in 2D.

\begin{Theorem}\label{thm:BKM criterion and GWP}
Let $c_0(x)\in H^s(\T^d)$ for $s>\f d 2+1$, and $(c,u)$ be a solution of (\ref{eq:HSCH}) stated in Theorem \ref{thm:local well-posedness}.
Let $T^*$ be the maximal existence time of the solution. If $T^*<+\infty$, then
\ben
\int_0^{T^*}\|\na c(t)\|_{L^\infty}^4dt=+\infty.
\een
In particular, this implies $T^*=+\infty$ for $d=2$. That is, the system (\ref{eq:HSCH}) is
globally well-posed in 2D.
\end{Theorem}

\no{\bf Proof.} First of all, we derive the basic energy law of the system.
Multiplying by $\mu$ on both sides of the third equation of (\ref{eq:HSCH}),
we get by integration by parts that
\beno
\int_{\T^d}c_t\mu dx+\int_{\T^d}u\cdot\na c \mu dx=-\f1{\textbf{Pe}}\int_{\T^d}|\na \mu|^2dx.
\eeno
Due to the definition of $\mu$, we have
\beno
\int_{\T^d}c_t\mu dx=\f d {dt}\Big(\int_{\T^d}f_0(c)dx+\f {\textbf{C}} 2\int_{\T^d}|\na c|^2dx\Big),
\eeno
and due to $\na\cdot u=0$,
\beno
\int_{\T^d}u\cdot\na c \mu dx=-\textbf{M}\int_{\T^d}u\cdot\big(\na p-\f1 {\textbf{M}}\mu\na c\big)dx
=12\textbf{M}\int_{\T^d}\eta(c)|u|^2dx.
\eeno
Thus we obtain the following classical energy equality \cite{LLG2002}
\beno
\f d {dt}\Big(\int_{\T^d}f_0(c)dx+\f {\textbf{C}} 2\int_{\T^d}|\na c|^2dx\Big)+
\f1{\textbf{Pe}}\int_{\T^d}|\na \mu|^2dx+12\textbf{M}\int_{\T^d}\eta(c)|u|^2dx=0.
\eeno
That is,
\ben\label{eq:energy equality}
E(t)+\f1{\textbf{Pe}}\int_0^t\int_{\T^d}|\na \mu(\tau)|^2dxd\tau
+12\textbf{M}\int_0^t\int_{\T^d}\eta(c)|u(\tau)|^2dxd\tau=E(0),
\een
where
\beno
E(t)\eqdefa \int_{\T^d}f_0(c(t,x))dx+\f {\textbf{C}} 2\int_{\T^d}|\na c(t,x)|^2dx.
\eeno
From the energy equality (\ref{eq:energy equality}), it follows that
\beno
\|c(t)\|_{H^1}^2+\f1{\textbf{Pe}}\int_0^t\|\na \mu\|_{L^2}^2d\tau\le E(0).
\eeno
On the other hand, we have
\beno
\|\na \Delta c\|_{L^2}\le C\big(\|\na \mu\|_{L^2}+\|\na c\|_{L^2}+\|c^2\na c\|_{L^2}\big).
\eeno
and by Sobolev inequality,
\beno
\|c^2\na c\|_{L^2}&\le& C\|c\|_{L^6}^2\|\na c\|_{L^6}\le C\|c\|_{H^1}^2\|c\|_{H^2}\\
&\le& C\|c\|_{H^1}^{\f52}\|c\|_{H^3}^\f12\le C\|c\|_{H^1}^{5}+
\f12\|c\|_{H^3},
\eeno
which implies that
\beno
\|c\|_{H^3}\le C\big(\|\na \mu\|_{L^2}+\|c\|_{H^1}+\|c\|_{H^1}^5\big).
\eeno
Therefore we conclude that
\ben\label{eq:H1-estimate}
\|c\|_{L^\infty(0,T; H^1)}+\|c\|_{L^2(0,T;H^3)}\le C\big(T,\|c_0\|_{H^1}\big).
\een

Next, we derive $H^2$ energy estimate of the solution. We have
\ben\label{eq:Energy-H2}
\f1 2\f {d}{dt}\|\Delta c\|_{L^2}^2+\f {\textbf{C}}{\textbf{Pe}}\|\Delta^2 c\|_{L^2}^2
&=&-\big(u\cdot\na c,\Delta^2 c\big)_{L^2}+\f1{\textbf{Pe}}\big(\Delta f_0'(c),\Delta^2c\big)\nonumber\\
&\le& \|u\|_{L^2}\|\na c\|_{L^\infty}\|\Delta^2 c\|_{L^2}
+\f1{\textbf{Pe}}\|\Delta f_0'(c)\|_{L^2}\|\Delta^2 c\|_{L^2}.
\een
It is easy to verify that
\beno
\|u\|_{L^2}&\le& C\big(\|\na p\|_{L^2}+\|\mu(c)\na c\|_{L^2}\big)\\
&\le& C\big(\|\na c\|_{L^\infty}\|\Delta c\|_{L^2}+(\|c\|_{L^3}+\|c\|_{L^9}^3)\|\na c\|_{L^6}\big)\\
&\le& C\big(\|\na c\|_{L^\infty}+\|c\|_{L^3}+\|c\|_{L^9}^3\big)\|c\|_{H^2},
\eeno
and
\beno
\|\Delta f_0'(c)\|_{L^2}\le C\big(1+\|c\|_{L^\infty}^2\big)\|c\|_{H^2}.
\eeno

Plugging them into (\ref{eq:Energy-H2}) yields that
\beno
\f {d}{dt}\|\Delta c\|_{L^2}^2+\|\Delta^2 c\|_{L^2}^2
\le C\big(1+\|\na c\|_{L^\infty}^4+\|c\|_{L^\infty}^4+\|c\|_{L^3}^4+\|c\|_{L^9}^{12}\big)
\|c\|_{H^2}^2,
\eeno
which along with Gronwall's inequality leads to
\ben\label{eq:Energy estimate-H2}
\|c\|_{H^2}\le \|c_0\|_{H^2}\exp\big(C\int_0^tH(\tau)d\tau\big),
\een
where
$
H(t)=1+\|\na c\|_{L^\infty}^4+\|c\|_{L^\infty}^4+\|c\|_{L^3}^4+\|c\|_{L^9}^{12}.
$

Now we are in position to prove the blow-up criterion. We will prove it by way of contradiction argument. Assume that
$T^*<+\infty$ and
\beno
\int_0^{T^*}\|\na c(t)\|_{L^\infty}^4dt<+\infty,
\eeno
which together with (\ref{eq:H1-estimate}) and Sobolev's inequality implies that
\beno
\int_0^{T^*}H(\tau)d\tau<+\infty,
\eeno
for example,
\beno
\int_0^{T^*}\|c(t)\|_{L^9}^{12}dt\le C\int_0^{T^*}\|c(t)\|_{H^1}^{11}\|c(t)\|_{H^3}dt<+\infty.
\eeno
Then we infer from (\ref{eq:Energy estimate-H2}) that
\beno
\|c\|_{L^\infty(0,T^*;H^2)}<+\infty,
\eeno
which implies that
\beno
\int_0^{T^*}G(t)dt<+\infty,\quad G(t)\quad \textrm{be as in Theorem }\ref{thm:local well-posedness}.
\eeno
Then the energy inequality (\ref{eq:energy estimate}) ensures that
\beno
\sup_{t\in [0,T^*]}\|c(t)\|_{H^s}^2+\int_0^{T^*}\|c(\tau)\|_{H^{s+2}}^2d\tau<+\infty,
\eeno
which means that the solution can be continued after $t=T^*$, and thus contradicts with the definition of
$T^*$.

As an application of blow-up criterion, we can deduce the global existence in 2D. Indeed, in two dimensional case,
we get by Gagliardo-Nirenberg inequality and (\ref{eq:H1-estimate}) that
\beno
\int_0^{T^*}\|\na c(t)\|_{L^\infty}^4dt\le C\int_0^{T^*}\|c(t)\|_{H^1}^2\|c(t)\|_{H^3}^2dt<+\infty,
\eeno
which implies $T^*=+\infty$ by the blow-up criterion.\ef

\section{A refined blow-up criterion in 3D}

 We first turn to a simple model relating to the Hele-Shaw-Cahn-Hilliard system:
\ben\label{eq:HSCH-simple}
\left\{
\begin{array}{l}
u=-\na p+\Delta c\na c,\quad \na\cdot u=0,\\
c_t+u\cdot\na c+\Delta^2 c=0.
\end{array}
\right.
\een
For this system, we still have the energy equality:
\beno
\|\na c(t)\|_{L^2}^2+2\int_0^t\|\na \Delta c(\tau)\|_{L^2}^2+\|u(\tau)\|^2_{L^2}d\tau=\|\na c_0\|_{L^2}.
\eeno
Moreover, if $c$ is a solution of (\ref{eq:HSCH-simple}), then
$c_\lambda(t,x)\eqdefa c(\lambda^4 t,\lambda x)$ is also a solution. It is easy to see that
\beno
\|\na c_\lambda(t,x)\|_{L^2}=\lambda^{\f d2-1}\|\na c(\lambda^4 t,x)\|_{L^2},\quad
\int_0^\infty\|\na \Delta c_\lambda(\tau)\|_{L^2}^2d\tau=\lambda^{2-d}\int_0^\infty\|\na \Delta c(\tau)\|_{L^2}^2d\tau.
\eeno
Thus, the energy is scaling invariance for $d=2$. From this view of point, the 2D system is critical
and the 3D system is supercritical like the 3D Navier-Stokes equations.
Due to the bi-Laplacian $\Delta^2$, there is no maximum principle for this system,
which is the main obstacle to obtain the global existence in 3D case. For the 2D critical QG equation
\beno
\theta_t+(-\Delta)^\f12\theta+u\cdot\na \theta=0,\quad u=\big(-(-\Delta)^{-\f12}\p_{x_2}\theta, (-\Delta)^{-\f12}\p_{x_1}\theta\big),
\eeno
Caffarelli and Vasseur \cite{Caf-Vas} proved the global regularity of weak solution.
The key step of their proof is to prove the H\"{o}lder continuity of the solution
by using the DeGiorgi method. Note that the QG equation has maximum principle.
For the 3D Hele-Shaw-Cahn-Hilliard system, we also show that the H\"{o}lder continuity of the solution will
control the blow-up of the solution.

\begin{theorem}\label{thm:refined blow-up criterion}
Let $\al \in(0,1)$ and $c_0(x)\in H^s(\T^3)$ for $s\ge 3$. Assume that $(c,u)$ be a solution of (\ref{eq:HSCH}) stated in Theorem \ref{thm:local well-posedness}.
Let $T^*$ be the maximal existence time of the solution. If $T^*<+\infty$, then
\beno
\int_0^{T^*}\|c(t)\|_{C^\al}^{\f 8 {\al} }dt=+\infty.
\eeno
\end{theorem}

\no{\bf Proof.}\,We will prove it by contradiction argument. Assume that
$T^*<+\infty$ and
\ben\label{eq:blow-up}
\int_0^{T^*}\|c(t)\|_{C^\al}^{\f 8 {\al} }dt<+\infty.
\een

Taking $\Delta_j$ to the third equation of (\ref{eq:HSCH}) to obtain
\beno
\p_t\Delta_jc+\f {\textbf{C}}{\textbf{Pe}}\Delta^2\Delta_jc
=-\Delta_j(u\cdot \na c)+\f 1 {\textbf{Pe}}\Delta \Delta_jf_0'(c).
\eeno
Making $L^2(\T^3)$ energy estimate, we get by Lemma \ref{lem:Berstein} that for $j\ge 0$,
\beno
\f{d}{dt}\|\Delta_j c\|_{L^2}^2+c2^{4j}\|\Delta_j c\|_{L^2}^2
\le C\big(\|\Delta_j(u\cdot \na c)\|_{L^2}+\|\Delta f_0'(c)\|_{L^2}\big)\|\Delta_j c\|_{L^2}.
\eeno
Dividing the above inequality by $\|\Delta_j c\|_{L^2}$ gives
\beno
\f{d}{dt}\|\Delta_j c\|_{L^2}+c2^{4j}\|\Delta_j c\|_{L^2}
\le C\big(\|\Delta_j(u\cdot \na c)\|_{L^2}+\|\Delta f_0'(c)\|_{L^2}\big),
\eeno
which implies that
\ben\label{eq:local L2}
\|\Delta_j c(t)\|_{L^2}\le \|\Delta_j c_0\|_{L^2}
+C\int_0^te^{-c2^{4j}(t-\tau)}\big(\|\Delta_j(u\cdot \na c)(\tau)\|_{L^2}+\|\Delta f_0'(c(\tau))\|_{L^2}\big)d\tau.
\een
We denote
\beno
\|c\|_{B^s_{2,\infty}}\eqdefa \sup_{j\ge -1}2^{js}\|\Delta_j c\|_{L^2}.
\eeno
Using the definition of Sobolev space, it is easy to find that
\beno
\|c\|_{H^{s-\epsilon}}^2\le \sum_{j\ge-1}2^{-2\varepsilon j}\|c\|_{B^s_{2,\infty}}^2\le C\|c\|_{B^s_{2,\infty}}^2
,\quad \forall \varepsilon>0.
\eeno

It follows from (\ref{eq:local L2}) that
\ben\label{eq:c-H3}
&&\|c(t)\|_{B^3_{2,\infty}}\le \|c(t)\|_{L^2}+\|c_0\|_{H^3}\nonumber\\
&&\qquad+C\sup_{j\ge 0}2^{3j}\int_0^te^{-c2^{4j}(t-\tau)}\big(\|\Delta_j(u\cdot \na c)(\tau)\|_{L^2}+\|\Delta f_0'(c(\tau))\|_{L^2}\big)d\tau.
\een
Now we claim that
\ben\label{eq:claim}
\|\Delta_j(u\cdot \na c)\|_{L^2}\le C2^{j(1-\al)}\|u\|_{L^2}\|c\|_{C^\al}.
\een
Now we have
\beno
\|u\|_{L^2}&\le& C\|\mu(c)\na c\|_{L^2}\le C\|c\|_{H^{3-\al}}\|c\|_{C^\al}+C\big(\|c\|_{L^3}+\|c\|_{L^6}^2\|c\|_{L^\infty}
\big)\|\na c\|_{L^6}\\
&\le& C\big(1+\|c\|_{H^1}+\|c\|_{H^1}^2\big)\|c\|_{C^\al}\|c\|_{B^{3}_{2,\infty}}.
\eeno
Here we used the product estimate
\beno
\|\Delta c \na c\|_{L^2}\le C\|c\|_{H^{3-\al}}\|c\|_{C^\al}\le C\|c\|_{B^{3}_{2,\infty}}\|c\|_{C^\al},
\eeno
which can be proved as in Lemma \ref{lem:product}. And similarly we have
\beno
\|\Delta f_0'(c)\|_{L^2}\le C\big(1+\|c\|_{C^\al}^2\big)\|c\|_{H^2}.
\eeno

Plugging the above estimates into (\ref{eq:c-H3}) yields that
\beno
&&\|c(t)\|_{B^3_{2,\infty}}\le \|c(t)\|_{L^2}+\|c_0\|_{H^3}\nonumber\\
&&\qquad+C\sup_{j\ge 0}2^{j(4-\al)}\int_0^te^{-c2^{4j}(t-\tau)}\big(1+\|c\|_{H^1}+\|c\|_{H^1}^2\big)
\big(1+\|c\|_{C^\al}^2\big)\|c\|_{B^3_{2,\infty}}d\tau,
\eeno
which along with H\"{o}lder inequality gives
\beno
&&\|c(t)\|_{L^\infty(0,t; B^3_{2,\infty})}\le \|c(t)\|_{L^\infty(0,t;L^2)}+\|c_0\|_{H^3}\\&&\qquad+
\big(1+\|c\|_{L^\infty(0,t;H^1)}+\|c\|_{L^\infty(0,t; H^1)}^2\big)
\big(t^\f 4 \al+\|c\|_{L^{\f 8\al}(0,t; C^\al)}^2\big)\|c\|_{L^\infty(0,t; B^3_{2,\infty})}.
\eeno
The above argument is still valid on the interval $[T,T^*)$ for $T<T^*$. Thus we get by
using (\ref{eq:H1-estimate}) that
\beno
&&\|c(t)\|_{L^\infty(T,T^*; B^3_{2,\infty})}\le \|c_0\|_{H^1}+\|c_0(T)\|_{H^3}\\&&\qquad+
C(\|c_0\|_{H^1})\big((T^*-T)^\f 4 \al+\|c\|_{L^{\f 8\al}(T,T^*; C^\al)}^2\big)\|c\|_{L^\infty(T,T^*; B^3_{2,\infty})}.
\eeno
Due to (\ref{eq:blow-up}), we can choose $T$ such that
\beno
C(\|c_0\|_{H^1})\big((T^*-T)^\f 4 \al+\|c\|_{L^{\f 8\al}(T,T^*; C^\al)}^2\big)\le \f12,
\eeno
Then we obtain
\beno
\|c(t)\|_{L^\infty(T,T^*; B^3_{2,\infty})}\le 2\big(\|c_0\|_{H^1}+\|c_0(T)\|_{H^3}\big),
\eeno
which implies by $\|\na c\|_{L^\infty}\le C\|c\|_{B^3_{2,\infty}}$ that
\beno
\int_0^{T^*}\|\na c(t)\|_{L^\infty}^4dt<+\infty,
\eeno
which is impossible by Theorem \ref{thm:BKM criterion and GWP} if $T^*<+\infty$.

It remains to prove (\ref{eq:claim}). As in proof of Lemma \ref{lem:product}, we have
\beno
\Delta_j(u\cdot \na c)&=&\Delta_j\sum_{|j-k|\le 4}S_{k-1}u\cdot\na \Delta_k c
+\Delta_j\sum_{|j-k|\le 4}\Delta_ku\cdot\na S_{k-1} c\\
&&+\Delta_j\sum_{|k-k'|\le 1,k\ge j-3}\Delta_k u\cdot\na\Delta_{k'}c=A_1+A_2+A_3.
\eeno
We get by Lemma \ref{lem:Berstein} that
\beno
\|A_1\|_{L^2}\le C\sum_{|j-k|\le 4}\|S_{k-1}u\|_{L^2}\|\na \Delta_k c\|_{L^\infty}
\le C2^{j(1-\al)}\|u\|_{L^2}\|c\|_{C^\al},
\eeno
and for $A_2$,
\beno
\|A_2\|_{L^2}&\le& C\sum_{|j-k|\le 4}\|\Delta_{k}u\|_{L^2}\|\na S_{k-1}c\|_{L^\infty}\\
&\le& C\|u\|_{L^2}\sum_{|j-k|\le 4}\sum_{\ell\le k-2}2^\ell\|\Delta_\ell c\|_{L^\infty}\\
&\le& C\|u\|_{L^2}\|c\|_{C^\al}\sum_{|j-k|\le 4}\sum_{\ell\le k-2}2^{\ell(1-\al)}
\le C2^{j(1-\al)}\|u\|_{L^2}\|c\|_{C^\al},
\eeno
and due to $\na \cdot u=0$,
\beno
\|A_3\|_{L^2}&\le& \|\Delta_j\sum_{|k-k'|\le 1,k\ge j-3}\na\cdot(\Delta_k u\Delta_{k'}c)\|_{L^2}\\
&\le& C2^j\sum_{|k-k'|\le 1,k\ge j-3}2^{-k'\al}\|u\|_{L^2}2^{k'\al}\|\Delta_{k'}c\|_{L^\infty}\\
&\le& C2^{j(1-\al)}\|u\|_{L^2}\|c\|_{C^\al}.
\eeno
Then the inequality (\ref{eq:claim}) follows from the estimates of $A_1, A_2$ and $A_3$. The proof of Theorem \ref{thm:refined blow-up criterion}
is completed. \ef

\section{Appendix} 
\setcounter{equation}{0}

Let us first recall some basic facts about the Littlewood-Paley theory.
Let $\varphi, \chi$ be two functions in $C^\infty(\T^d)$ such that
$\textrm{supp}\widehat{\varphi} \subset \{\frac{3}{4}\le|\xi|\le\frac{8}{3}\}$,
$\textrm{supp}\widehat{\chi}\subset\{|\xi|\le\frac{4}{3}\}$ and
\begin{align*}
\widehat\chi(\xi)+\sum_{j\ge 0}\widehat\varphi(2^{-j}\xi)=1.
\end{align*}
Then the Littlewood-Paley operators are defined by
\beno
&&\Delta_jf=\varphi_j\ast f=\int_{\T^d}\varphi_j(x-y)f(y)dy,\quad \varphi_j(x)=2^{jd}\varphi(2^jx),\quad j\ge 0,\\
&&S_jf=\chi_j\ast f=\sum_{k=-1}^{j-1}\Delta_kf,\quad \Delta_{-1}f=\chi\ast f.
\eeno
Some classical spaces can be characterized in terms of $\Delta_j$.
Let $s\in \R$, the Sobolev space $H^s({\T^d})$ is defined by
\beno
H^s(\T^d)\eqdefa \big\{u\in \cD'(\T^d):
\|u\|_{H^s}^2\eqdefa \sum_{j\ge -1}2^{2js}\|\Delta_j u\|_{L^2}^2<\infty\big\}.
\eeno
We denote by $(u,v)_{H^s}$ the inner product in $H^s(\T^d)$.
And for $s\in (0,1)$, the H\"{o}lder space $C^s(\T^d)$ is defined by
\beno
C^s(\T^d)\eqdefa \big\{u\in \cD'(\T^d):
\|u\|_{C^s}\eqdefa \sup_{j\ge -1}2^{js}\|\Delta_j u\|_{L^\infty}\big\}.
\eeno
We refer to \cite{Triebel} for more details.  Let us recall Bony's decomposition from \cite{Bony}:
\ben\label{Bony}
fg=T_fg+T_gf+R(f,g),
\een
where
\beno
T_fg=\sum_{j\ge -1} S_{j-1}f\Delta_jg,\quad R(f,g)=\sum_{|j-j'|\le 1}\Delta_{j}f\Delta_{j'}g.
\eeno
We also denote $\widetilde{R}(f,g)=T_gf+R(f,g)$.

\begin{lemma}\label{lem:Berstein}\cite{Chemin}
Let $k\in \N, 1\le p\le q\le\infty$. Then there exists a positive constant $C$ independent of $j$ such that
\beno
&&\|\partial^\alpha
\Delta_j f\|_{L^q}+\|\partial^\alpha
S_j f\|_{L^q}\le
C2^{j{|\alpha|}+dj(\frac{1}{p}-\frac{1}{q})}\|f\|_{L^p},
\\
&&\|\Delta_j f\|_{L^p}\le C2^{-jk}\sup_{|\al|=k}\|\partial^\al \Delta_j f\|_{L^p},\quad j\ge 0. \eeno
\end{lemma}

\begin{lemma}\label{lem:product}
Let $s\ge 0$. Then there holds
\ben\label{eq:product}
\|fg\|_{H^s}\le C\big(\|f\|_{L^\infty}\|g\|_{H^s}+\|f\|_{H^s}\|g\|_{L^\infty}\big).
\een
If $0<\sigma\le \f d 2$, then there holds
\ben\label{eq:product-refine}
\|fg\|_{H^s}\le C\big(\|f\|_{H^s}\|g\|_{L^\infty}+\|f\|_{H^{\f d2-\sigma}}\|g\|_{H^{s+\sigma}}\big).
\een
\end{lemma}

\no{\bf Proof.} The inequality  (\ref{eq:product}) is classical, see \cite{Kato}. Here we only present the proof of (\ref{eq:product-refine}).
Using the Bony's decomposition (\ref{Bony}) to write
\beno
\Delta_j(fg)=\Delta_j(T_{f}g)+\Delta_j(T_{g}f)+\Delta_jR(f,g).
\eeno
Taking into considering the support of Fourier transform of the term $T_{f}g$,
we have
\beno
\Delta_j(T_{f}g)=\sum_{|j'-j|\le4}\Delta_j(S_{j'-1}f\Delta_{j'}g).
\eeno
Due to $0<\sigma\le \f d2$, this gives by Lemma \ref{lem:Berstein} that
\beno
\|S_{j}f\|_{L^\infty}\le \left\{\begin{array}{l}
C2^{j\f d2}\|f\|_{L^2},\quad \textrm{if } \sigma=\f d2,\\
C\displaystyle\sum_{k\le j-1}2^{k\f d2}\|\Delta_k f\|_{L^2}\le C2^{j\sigma}\|f\|_{H^{\f d2-\sigma}}, \quad \textrm{if }\sigma<\f d2,
\end{array}\right.
\eeno
which implies that
\ben\label{eq:LH-est} \|\Delta_j(T_{f}g)\|_{L^2}&\le& C
\sum_{|j'-j|\le4}\|S_{j'-1}f\|_{L^\infty}
\|\Delta_{j'}g\|_{L^2}\nonumber\\
&\le& C\|f\|_{H^{\f d2-\sigma}}\sum_{|j'-j|\le 4}2^{j'\sigma}\|\Delta_{j'}g\|_{L^2}\nonumber\\
&\le& C2^{-js}c_j\|f\|_{H^{\f d2-\sigma}}\|g\|_{H^{s+\sigma}},
\een
here and hereafter $\{c_j\}_{j\ge -1}$ denotes a sequence satisfying $\|\{c_j\}_{j\ge -1}\|_{\ell^2}\le 1$.

Similarly, we have
\ben\label{eq:HL-est}
\|\Delta_j(T_{g}f)\|_{L^2}&\le& C
\sum_{|j'-j|\le4}\|S_{j'-1}g\|_{L^\infty}
\|\Delta_{j'}f\|_{L^2}\nonumber\\
&\le& C\sum_{|j'-j|\le4}\|g\|_{L^\infty}\|\Delta_{j'}f\|_{L^2}\nonumber\\
&\le& C2^{-js}c_j\|g\|_{L^\infty}\|f\|_{H^{s}}.
\een
Noticing that, after taking into account the support of the Fourier transforms,
\beno
\Delta_jR(f,g)=\sum_{j',j''\ge
j-3;|j'-j''|\le 1}\Delta_j(\Delta_{j'}f\Delta_{j''}g),
\eeno
it follows from Lemma \ref{lem:Berstein} that
\ben\label{eq:HH-est}
\|\Delta_jR(f,g)\|_{L^2} &\le& C
\sum_{j',j''\ge j-3;|j'-j''|\le 1}2^{j\f d2}\|\Delta_{j'}f\|_{L^2}\|\Delta_{j''}g\|_{L^2}\nonumber\\
&\le&  C2^{-js}\sum_{j',j''\ge j-3;|j'-j''|\le 1}2^{(j-j')(\f d2+s)}2^{j'(\f d2-\sigma)}\|\Delta_{j'}f\|_{L^2}
2^{j''(s+\sigma)}\|\Delta_{j''}g\|_{L^2}\nonumber\\
&\le& C2^{-js}c_j\|f\|_{H^{\f d2-\sigma}}\|g\|_{H^{s+\sigma}}.
\een
Thanks to the definition of Sobolev space, (\ref{eq:product-refine}) follows from (\ref{eq:LH-est})-(\ref{eq:HH-est}).\ef

\begin{lemma}\label{lem:nonlinear}\cite{Triebel}
Let $s>0$. Assume that $F(\cdot)$ is a smooth function on $\R$ with $F(0)=0$. Then we have
\beno
\|F(f)\|_{H^s}\le C(1+\|f\|_{L^\infty})^{\lfloor s \rfloor+1}\|f\|_{H^s},
\eeno
where the constant $C$ depends on $\displaystyle\sup_{k\le \lfloor s\rfloor+2, |t|\le \|f\|_{L^\infty}}\|F^{(k)}(t)\|_{L^\infty}$.
\end{lemma}

\begin{lemma}\label{lem:commutator}
Let $s>0$. Then there holds
\beno
\|\langle D\rangle^s(fg)-f\langle D\rangle^sg\|_{L^2}\le C\big(\|f\|_{H^{s+2}}\|g\|_{L^2}+\|f\|_{H^2}\|g\|_{H^{s-\f12}}\big).
\eeno
If $s\in (0,1]$, then we have
\beno
\|\langle D\rangle^s(fg)-f\langle D\rangle^sg\|_{L^2}\le C\|f\|_{H^{s+2}}\|g\|_{L^2}.
\eeno
Here the Fourier multiplier $\langle D\rangle^s$ is defined by
\beno
\langle D\rangle^s f(x)=\sum_{k\in \Z^d}(1+|k|^2)^\f s2e^{2\pi ik\cdot x}\widehat f(k).
\eeno
\end{lemma}

\no{\bf Proof.} Using Bony's decomposition (\ref{Bony}) to write
\beno
&&\langle D\rangle^s(fg)=\langle D\rangle^s(T_fg)+\langle D\rangle^sT_gf+\langle D\rangle^sR(f,g),\\
&&f\langle D\rangle^sg=T_f\langle D\rangle^sg+T_{\langle D\rangle^sg}f+R(f,\langle D\rangle^sg).
\eeno
Thus we have
\beno
\langle D\rangle^s(fg)-f\langle D\rangle^sg=\langle D\rangle^s(T_fg)-T_f\langle D\rangle^sg+\pi(f,g),
\eeno
where
\beno
\pi(f,g)=\langle D\rangle^sT_gf+\langle D\rangle^sR(f,g)-T_{\langle D\rangle^sg}f-R(f,\langle D\rangle^sg).
\eeno

As in the proof of (\ref{eq:product-refine}), we can deduce by Lemma \ref{lem:Berstein} that
\beno
\|\pi(f,g)\|_{L^2}\le C\|f\|_{H^{s+2}}\|g\|_{L^2}.
\eeno
We illustrate the process by working out the estimate on the first term. Thanks to Lemma \ref{lem:Berstein}, we have
\beno
\|\langle D\rangle^sT_gf\|_{L^2}^2&=&\sum_{j\ge -1}\|\Delta_j \langle D\rangle^sT_gf\|_{L^2}^2
\le C\sum_{j\ge -1}2^{2js}\|\Delta_jT_gf\|_{L^2}^2\\
&\le& C\sum_{|j-j'|\le 4}2^{2js}\|S_{j'-1}g\Delta_{j'}f\|_{L^2}^2\\
&\le& C\sum_{|j-j'|\le 4}2^{2js}\|S_{j'-1}g\|_{L^\infty}^2\|\Delta_{j'}f\|_{L^2}^2\\
&\le& C\sum_{|j-j'|\le 4}2^{2j(s+\f d2)}\|g\|_{L^2}^2\|\Delta_{j'}f\|_{L^2}^2\\
&\le& C\|g\|_{L^2}^2\|f\|_{H^{s+\f d2}}^2\le C\|g\|_{L^2}^2\|f\|_{H^{s+2}}^2.
\eeno

Let $m(\xi_1,\xi_2)$ be the symbol of the paraproduct operator $T_fg$.
Then $\langle D\rangle^s(T_fg)-T_f\langle D\rangle^sg$ has the symbol
\beno
m(\xi_1,\xi_2)\big(\langle \xi_1+\xi_2\rangle^s-\langle\xi_2\rangle^s\big),
\eeno
which is supported in the region $|\xi_1+\xi_2|\sim |\xi_2|$. By the fundamental theorem of calculus we have
\beno
m(\xi_1,\xi_2)\big(\langle \xi_1+\xi_2\rangle^s-\langle\xi_2\rangle^s\big)
=\int_0^1\xi_1\cdot m(\xi_1,\xi_2)\na h^s(t\xi_1+\xi_2)dt,\quad h^s(\xi)=\langle \xi\rangle^s.
\eeno
It is easy to verify that
$\langle\xi_1\rangle^{\theta}m(\xi_1,\xi_2)\na h^s(t \xi_1+\xi_2)\langle\xi_2\rangle^{1-\theta-s}$ with
$\theta\in [0,1]$ is a Coifman-Meyer paraproduct uniformly for $t\in [0,1]$. Then we have
\beno
\|\langle D\rangle^s(T_fg)-T_f\langle D\rangle^sg\|_{L^2}\le
C\|\langle D\rangle^{1-\theta}f\|_{L^p}\|\langle D\rangle^{s+\theta-1}g\|_{L^q}
\eeno
for $\theta\in [0,1]$, $\f1p+\f1q=\f12$ and $1<q<\infty$, see P. 106 in \cite{Work}. Taking $\theta=\f12$ ,$(p,q)=(\infty,2)$ for $d=2$, and
$\theta=0$, $(p,q)=(6,3)$ for $d=3$, we obtain
\beno
\|\langle D\rangle^s(T_fg)-T_f\langle D\rangle^sg\|_{L^2}\le C\|f\|_{H^2}\|g\|_{H^{s-\f12}}.
\eeno
In case of $s\in (0,1]$, taking $\theta=1-s$ and $(p,q)=(\infty,2)$ to obtain
\beno
\|\langle D\rangle^s(T_fg)-T_f\langle D\rangle^sg\|_{L^2}\le C\|f\|_{H^{s+2}}\|g\|_{L^2}.
\eeno
This completes the proof of Lemma \ref{lem:commutator}.\ef

\bigskip

\noindent {\bf Acknowledgments.} Part of the work was carried out while the authors were
long term visitors at IMA at University of Minnesota. The hospitality and support of IMA are graciously acknowledged.
Xiaoming Wang is supported in part by NSF. He also acknowledges helpful conversation with David Ambrose, Maurizio Grasselli, Xiaoqiang Wang and Steve Wise.
Zhifei Zhang is
supported by NSF of China under Grant 10990013 and 11071007.


\end{document}